\documentclass{article}

\usepackage{amsfonts,latexsym}
\usepackage{graphicx,graphics,hieroglf,skak}
\usepackage{amsthm}
\usepackage{amssymb, amsmath}
\usepackage{url,enumitem}

\newtheorem{theorem}{Theorem}
\newtheorem{property}[theorem]{Property}

\newtheorem*{problem}{Problem}

\begin{document}

\title{The parbelos, a parabolic analog of the arbelos}
\author{Jonathan Sondow}
\date{}
\maketitle

\begin{abstract}
The \emph{arbelos} is a classical geometric shape bounded by three mutually tangent semicircles with collinear diameters. We introduce a parabolic analog, the \emph{parbelos}. After a review of the parabola, we use theorems of Archimedes and Lambert to demonstrate seven properties of the parbelos, drawing analogies to similar properties of the arbelos, some of which may also be new.
\end{abstract}

\section{Introduction: The Arbelos and the Parbelos.}

\begin{figure}[ht]
\begin{minipage}[b]{0.4\linewidth}
\centering
\includegraphics[width=2.2in]{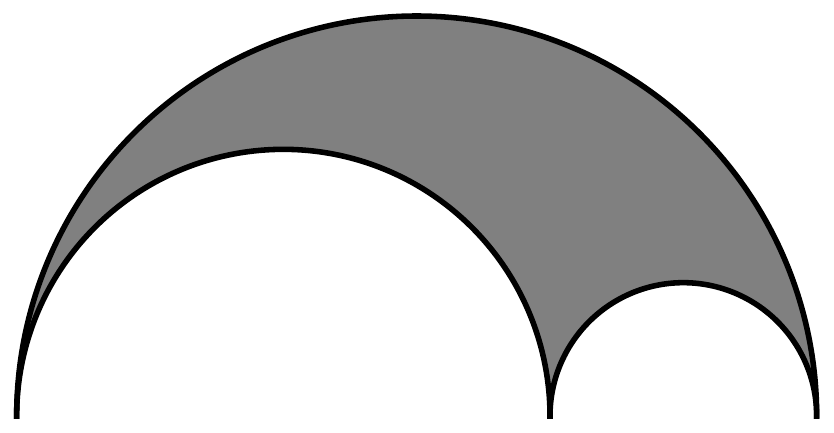}
\caption{The arbelos.}
\label{FIG:arbelos}
\end{minipage}
\hspace{1cm}
\begin{minipage}[b]{0.4\linewidth}
\centering
\includegraphics[width=2.2in]{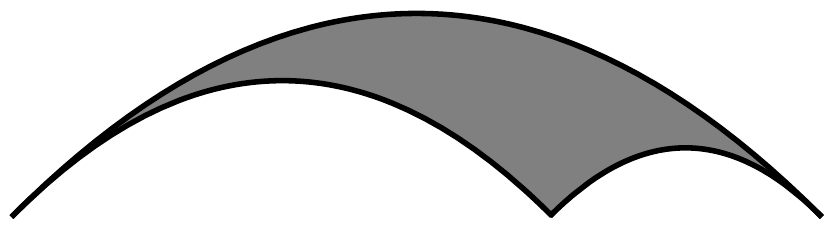}
\caption{The parbelos.}
\label{FIG:parbelos}
\end{minipage}
\end{figure}

The \emph{arbelos} or \emph{shoemaker's knife} is a classic figure from Greek geometry bounded by three pairwise tangent semicircles with diameters lying on the same line. (See Figure~\ref{FIG:arbelos}.)
There is a long list of remarkable properties of the arbelos---consult Boas's survey \cite{boas}, its forty-four references,
and Bogomolny's website \cite{bogomolnyA}.

Just as all circles are similar, so too all parabolas are similar. (See, e.g., \cite[p.~118]{bls}.
The same is not true for the other conic sections, because the similarity class of an ellipse or a hyperbola depends on its eccentricity.) For that reason,
one might expect to find a parabolic analog of the arbelos in the literature of the past two millennia. However, extensive searches have 
failed to uncover any mention of one.
This note provides and studies such an analog.

To define it, recall first that the latus rectum of a parabola is the focal chord parallel to the directrix (see Section~2). Now, replace the semicircles of the arbelos with the latus rectum arcs of parabolas, all opening in the same direction, whose foci are the centers of the semicircles.
The region bounded by the three arcs is the \emph{parbelos} associated to the arbelos. (See Figure~\ref{FIG:parbelos}.)

\noindent

More intrinsically, given three points $\mathsf{C_1,C_2,C_3}$ on a line, construct parabolas that open in the same direction and have latera recta $\mathsf{C_1C_2,C_2C_3,C_1C_3}$. The three latus rectum arcs enclose the parbelos. The points $\mathsf{C_1,C_2,C_3}$ are its \emph{cusps}.

Unlike in the arbelos, the arcs of the parbelos are not pairwise tangent: the inner two are tangent to the outer one, but not to each~other, as we show.

Section~2 is a review of the parabola. In Section~3 we use theorems of Archimedes and Lambert to prove seven properties of the parbelos, drawing analogies to similar properties found in the arbelos, some of which may also be new. The seventh constructs a parbelos directly from an arbelos via a locus.
Along the way we mention the Universal Parabolic Constant (an analog of $\pi$), Newton's teacher Barrow, and an origami fold.

\section{The Parabola.}

Recall that a parabola~$\mathsf{P}$ is the locus of points equidistant from a point~$\mathsf{F}$, called the \emph{focus}, and a line~$\mathsf{L}$, the \emph{directrix}. The distance $p>0$ from~$\mathsf{F}$ to~$\mathsf{L}$ is the \emph{focal parameter}. The point~$\mathsf{V}$ at a distance $a:=\frac p2$ from both $\mathsf{F}$ and~$\mathsf{L}$ is the \emph{vertex} of~$\mathsf{P}$. The chord $\mathsf{C_1C_2}$ of~$\mathsf{P}$ passing through~$\mathsf{F}$ and parallel to~$\mathsf{L}$ is the \emph{latus rectum}. Since $\mathsf{C_1}$ and $\mathsf{C_2}$ lie on~$\mathsf{P}$, the length of $\mathsf{C_1C_2}$ equals $2p=4a$; one half of that is the \emph{semi-latus rectum} $p=2a$. The arc of~$\mathsf{P}$ with endpoints $\mathsf{C_1}$ and $\mathsf{C_2}$ is the \emph{latus rectum arc}.

\begin{figure}[htbp]
	\begin{center}
		\includegraphics[width=2.45in]{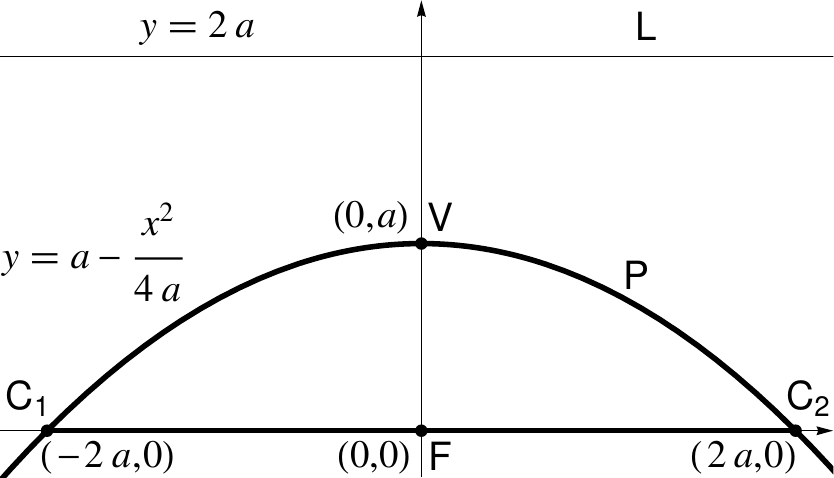}
		\caption{The downward-opening parabola with latus rectum $[-2a,2a]$.}
		\label{FIG:parabola}
	\end{center}
\end{figure}

These notations are illustrated in Figure~\ref{FIG:parabola}, which shows the unique downward-opening parabola~$\mathsf{P}$ whose latus rectum is the interval $[-2a,2a]$ on the \mbox{$x$-axis}. 
The equation of~$\mathsf{P}$ is $y=a-\frac{x^2}{4a}$.

Just as the ratio of the length of any semicircle to its radius is always~$\pi$, the ratio of the length~$s$ of the latus rectum arc of any parabola to its semi-latus rectum (and focal parameter)~$p$
is also a constant, namely,
the \emph{Universal Parabolic Constant} $P=\frac{s}{p}$---see Reese and Sondow \cite{rs}.
For a geometric method of computing the length of a parabolic arc, due to Isaac Barrow (1630--1677), see D\"{o}rrie \cite[Problem 58]{dorrie}; for the ``sweeping tangents'' method, see Apostol and Mnatsakanian \cite{am}.

\section{Properties of the Parbelos.}

Returning to the parbelos, we describe seven properties of it. The first two are analogous to classical ones of the arbelos, and are based on similarity.

\begin{property} \label{PROP:boundaries}
The upper and lower boundaries of the parbelos have the same length. 
\end{property}

\begin{proof}
(We paraphrase Boas \cite[p.~237]{boas} on the arbelos.) This is immediate from the knowledge that the length of the latus rectum arc of a parabola is proportional to its semi-latus rectum; one does not even need to know that the constant of proportionality is the Universal Parabolic Constant~$P$.
\end{proof}

Our second property of the parbelos is the direct analog of part of a deeper property of the arbelos discovered by Schoch in 1998---see \cite[Figure~6]{schoch}.

\begin{figure}[htbp]
	\begin{center}
		\includegraphics[width=2.4in]{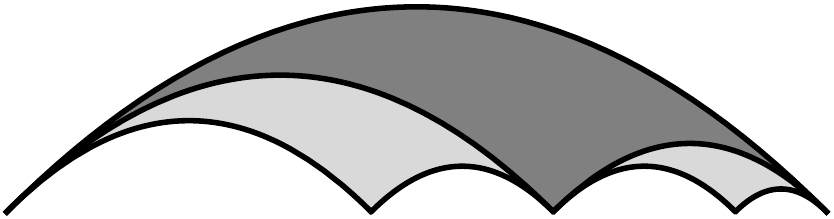}
		\caption{Three similar parbeloses.}
		\label{FIG:SimilarParbeloses}
	\end{center}
\end{figure}

\begin{property} \label{PROP:similar}
Under each lower arc of the parbelos, construct a new parbelos similar to the original $($see Figure~\ref{FIG:SimilarParbeloses}$)$. Of the four new lower arcs, the middle two are congruent, and their common length equals one half the harmonic mean of the lengths of the original lower arcs.
\end{property}

\begin{proof} Denote by $\ell_L$ and $\ell_R$ the lengths of the original left and right lower arcs, and by $\ell_1,\ell_2,\ell_3,\ell_4$ the lengths of the four new lower arcs. 
By similarity, we have the equalities
$$\ell_2=\frac{\ell_L}{\ell_L+\ell_R}\cdot \ell_R =\frac12\cdot\frac{2}{\frac{1}{\ell_L}+\frac{1}{\ell_R}} =\frac{\ell_R}{\ell_L+\ell_R}\cdot \ell_L= \ell_3.$$
As the arcs are latus rectum arcs, the result follows.
\end{proof}

\begin{figure}[ht]
\begin{minipage}[b]{0.4\linewidth}
\centering
\includegraphics[width=2.3in]{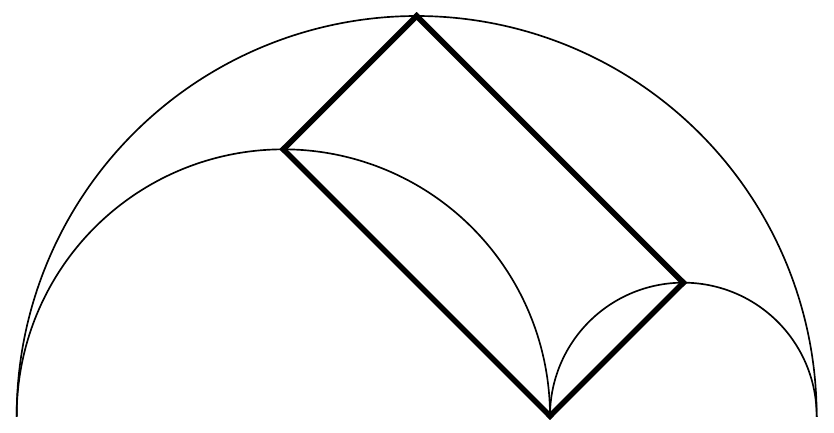}
\caption{The cusp-midpoints rectangle of the arbelos.}
\label{FIG:arbelosrectangle}
\end{minipage}
\hspace{1cm}
\begin{minipage}[b]{0.4\linewidth}
\centering
\includegraphics[width=2.3in]{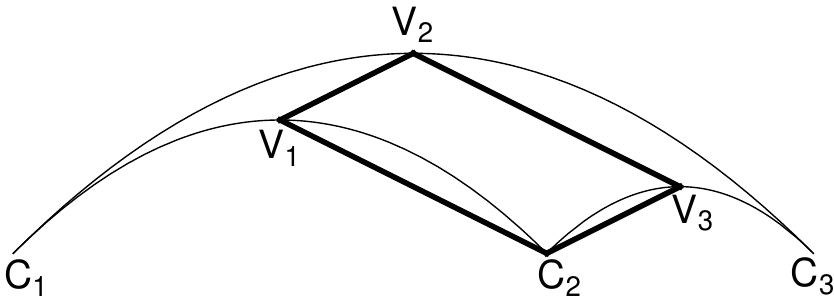}
\caption{The cusp-vertices parallelogram of the parbelos.}
\label{FIG:parallelogram}
\end{minipage}
\end{figure}

The next property of the parbelos is analogous to the fact that \emph{the area of the arbelos equals $\pi/2$ times the area of its cusp-midpoints rectangle}, determined by the middle cusp and the midpoints of the three semicircular arcs (see Figure~\ref{FIG:arbelosrectangle}).
The proof of this fact is similar to that of Property~\ref{PROP:parallelogram}.

\begin{property} \label{PROP:parallelogram}
The middle cusp of the parbelos and the vertices of its three parabolas determine a parallelogram, the {\em cusp-vertices parallelogram}. The area of the parbelos equals $4/3$ times the area of its cusp-vertices parallelogram.
\end{property}

\noindent

\begin{proof}
A glance at Figure~\ref{FIG:parabola} reveals that, in Figure~\ref{FIG:parallelogram}, the cusp $\mathsf{C_1}$ and the vertices $\mathsf{V_1}$ and $\mathsf{V_2}$ all lie on a line with slope $1/2$, as do the cusp $\mathsf{C_2}$ and the vertex $\mathsf{V_3}$. 
Likewise, $\mathsf{C_3}, \mathsf{V_3},$ and $\mathsf{V_2}$ all lie on a line with slope $-1/2$, as do $\mathsf{C_2}$ and $\mathsf{V_1}$. The first statement follows.

Let us denote the area of a parabolic segment by~$\betteris$, and that of a triangle by~$\triangle$. The area of the parbelos is then
\begin{align*}
\text{area}\,\mathsf{C_1V_2C_3V_3C_2V_1} = \betteris \mathsf{C_1V_2C_3} - (\betteris \mathsf{C_1V_1C_2} + \betteris \mathsf{C_2V_3C_3})
\end{align*}
and the area of the cusp-vertices parallelogram is
\begin{align*}
\text{area}\,\mathsf{C_2V_1V_2V_3} = \triangle \mathsf{C_1V_2C_3} - (\triangle \mathsf{C_1V_1C_2} + \triangle \mathsf{C_2V_3C_3} ).
\end{align*}
By Archimedes's calculation of the area of a parabolic segment (or by integral calculus), each~$\betteris$ equals 4/3 times the corresponding~$\triangle$. This proves the second statement.
\end{proof}

For an exposition of Archimedes of Syracuse's quadrature of the parabola (nineteen centuries before the Fundamental Theorem of Calculus!), see \cite[Problem~56]{dorrie}.

Here is another characterization of the area of the parbelos.

\begin{property} \label{PROP:rectangle}
The four tangents to the parbelos at its three cusps enclose a rectangle, the \emph{tangent rectangle}. The parbelos has two thirds the area of its tangent rectangle.
\end{property}

\begin{figure}[ht]
\begin{minipage}[b]{0.45\linewidth}
\centering
\includegraphics[width=2.3in]{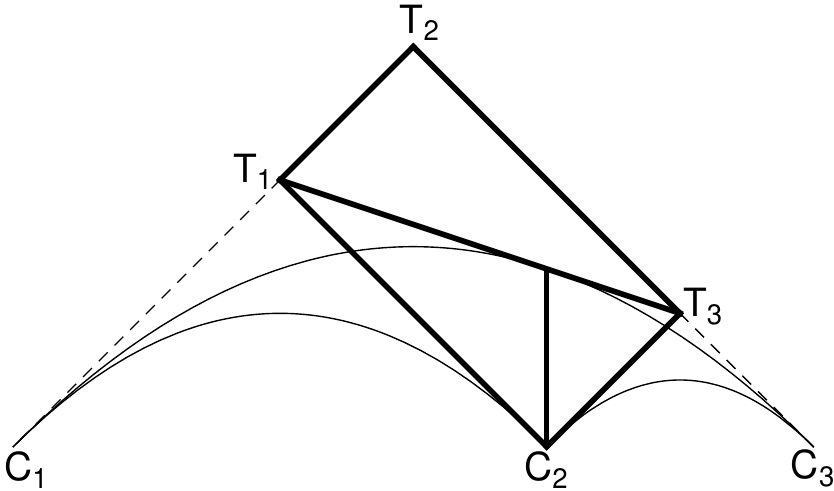}
\caption{The tangent rectangle of the parbelos, a diagonal, and an angle bisector.}
\label{FIG:RectangleAll}
\end{minipage}
\hspace{.1cm}
\begin{minipage}[b]{0.45\linewidth}
\centering
\includegraphics[width=2.3in]{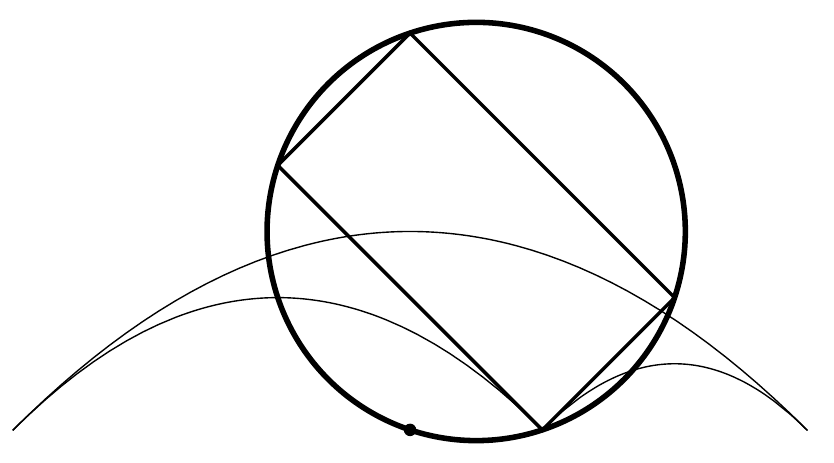}
\caption{The circumcircle of the tangent rectangle and the focus of the upper parabola.}
\label{FIG:RectangleCircle}
\end{minipage}
\end{figure}

\begin{proof}
By Figure~\ref{FIG:parabola} and calculus (or by the parabola's reflection property), the latus rectum of a parabola makes an angle of $\pi/4$ with the tangent line at each endpoint.
The first statement follows.

The area of the tangent rectangle in Figure~\ref{FIG:RectangleAll} is
\begin{align*}
\text{area}\, \mathsf{C_2T_1T_2T_3} = \triangle \mathsf{C_1T_2C_3} - (\triangle \mathsf{C_1T_1C_2} + \triangle \mathsf{C_2T_3C_3} ).
\end{align*}
From Figure~\ref{FIG:parabola}, each~$\triangle$ equals twice the corresponding~$\triangle$ in the previous proof, and we are done.
\end{proof}

The tangent rectangle also figures in the following two properties.

\begin{property} \label{PROP:diagonal}
In the tangent rectangle of the parbelos, the diagonal opposite the cusp is tangent to the upper parabola. The contact point lies on the bisector of the angle at the cusp.
\end{property}

\begin{proof}
In Figure~\ref{FIG:RectangleAll}, choose coordinates $\mathsf{C_1}=(0,0)$, $\mathsf{C_2}=(2b,0)$, and $\mathsf{C_3}=(4a,0)$.
The endpoints of the diagonal opposite $\mathsf{C_2}$ are then $\mathsf{T_1}=(b,b)$ and $\mathsf{T_3}=(2a+b,2a-b)$. Now, the equation of the line through $\mathsf{T_1}$ and $\mathsf{T_3}$, and the equation of the upper parabola, are
\begin{equation*} \label{EQ:tan}
	y =f(x):=\frac{a-b}{a}x + \frac{b^2}{a} \qquad \text{and} \qquad  y =g(x):=a-\frac{(x-2a)^2}{4a}.
\end{equation*}
Setting $f(x)=g(x)$, the only solution is $x=2b$. As $g'(2b)=\frac{a-b}{a}$, we infer Property~\ref{PROP:diagonal}.
\end{proof}

\begin{problem}
Find a proof by synthetic or Euclidean geometry, without introducing Cartesian coordinates.
\end{problem}
\noindent(\emph{Added in proof:} Emmanuel Tsukerman \cite{tsukerman}, an undergraduate (!) at Stanford University, has solved the Problem. He gives a synthetic proof of Property~\ref{PROP:diagonal} using a converse---which he proves---to Lambert's Theorem below.)

\smallskip
Property~\ref{PROP:diagonal} has a surprising consequence.

\begin{property} \label{PROP:circumcircle}
The circumcircle of the tangent rectangle of the parbelos passes through the focus of the upper parabola $($see Figure~\ref{FIG:RectangleCircle}$)$.
\end{property}

We give two proofs. The first uses the statement of Property~\ref{PROP:diagonal}; the second uses its proof.

\begin{proof}[Proof 1.]
Any three lines tangent to a parabola bound a \emph{tangent triangle}.
By Property~\ref{PROP:diagonal}, a diagonal and two adjacent sides of the tangent rectangle form a tangent triangle 
of the upper parabola. 
Property~\ref{PROP:circumcircle} now follows from \emph{Lambert's Theorem}, which asserts that \emph{the circumcircle of any tangent triangle of a parabola passes through the focus.} 
\end{proof}

\begin{proof}[Proof 2.]
The center of the circumcircle is $\frac12(\mathsf{T_1}+\mathsf{T_3})=(a+b,a)$, which is equidistant from the cusp $\mathsf{C_2}=(2b,0)$ and the focus $\mathsf{F}=(2a,0)$ of the upper parabola. This proves the property.
\end{proof}

See \cite[Problem~44]{dorrie} or \cite{tsukerman} for a proof of Lambert's Theorem. (Johann Heinrich Lambert (1728--1777) was a Swiss mathematician, physicist, and astronomer who gave the first proof that $\pi$~is irrational.)

\begin{figure}[ht]
\begin{minipage}[b]{0.45\linewidth}
\centering
\includegraphics[width=2.2in]{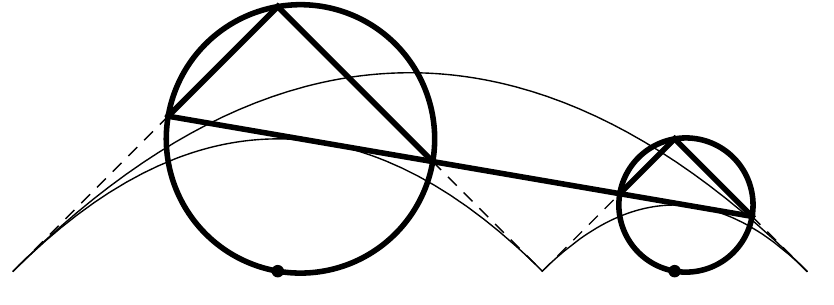}
\caption{The common tangent to the lower parabolas of the parbelos, their foci, and the circumcircles of similar tangent triangles.}
\label{FIG:TwoCircumcircles}
\end{minipage}
\hspace{.1cm}
\begin{minipage}[b]{0.45\linewidth}
\centering
\includegraphics[width=2.2in]{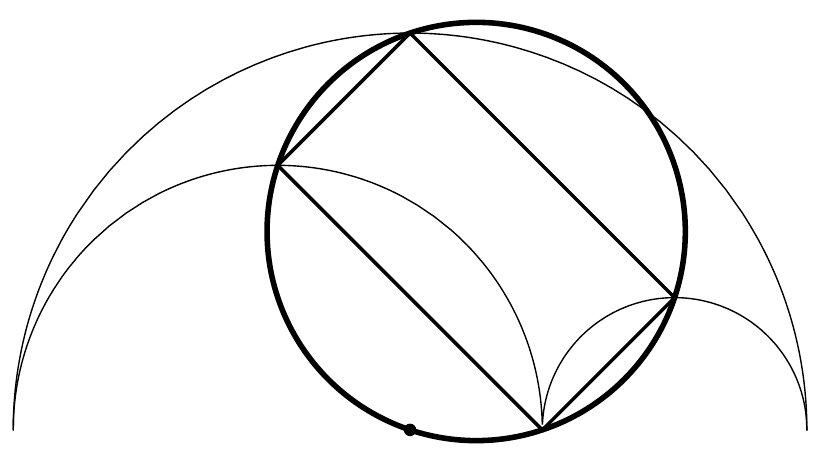}
\caption{The circumcircle of the cusp-midpoints rectangle of the arbelos and the center of the upper semicircle.}
\label{FIG:ArbelosParbelos}
\end{minipage}
\end{figure}

Of course, his theorem can also be applied to tangent triangles of the \emph{lower} arcs in the parbelos. For example, Figure~\ref{FIG:TwoCircumcircles} shows two similar tangent triangles formed by the tangents at the cusps and a line constructed tangent to both lower parabolas. The construction requires either solving a cubic equation or doing a Beloch origami fold---see Hull~\cite{hull}.

\noindent

Using Figures~\ref{FIG:arbelosrectangle} and~\ref{FIG:RectangleCircle}, it is easy to show that \emph{the cusp-midpoints rectangle of the arbelos coincides with the tangent rectangle of the associated parbelos}, whose foci are the centers of the semicircles of the arbelos. Hence, by Property~\ref{PROP:circumcircle},
\emph{the circumcircle of the cusp-midpoints rectangle of the arbelos passes through the center of the upper semicircle}---see Figure~\ref{FIG:ArbelosParbelos}.

Our final property describes how to construct a parbelos directly from an arbelos via a locus.

\begin{property} \label{PROP:locus}
The locus of the centers of circles inscribed in a semicircle of the arbelos is the boundary of a parbelos with its cusps deleted. The arbelos and parbelos share the same cusps.
\end{property}

\begin{proof}
We claim that \emph{the locus of the centers of circles inscribed in any semicircle 
is the open latus rectum arc of a parabola whose latus rectum is the diameter of the semicircle}. A more general fact was discovered by Byer, Lazebnik, and Smeltzer \cite[p.~118]{bls}. (\emph{Added in proof:} Tahir \cite[p.~30]{tahir} discovered it eleven years earlier.) We adapt their elegant proof as follows.

\begin{figure}[htbp]
	\begin{center}
		\includegraphics[width=2.9in]{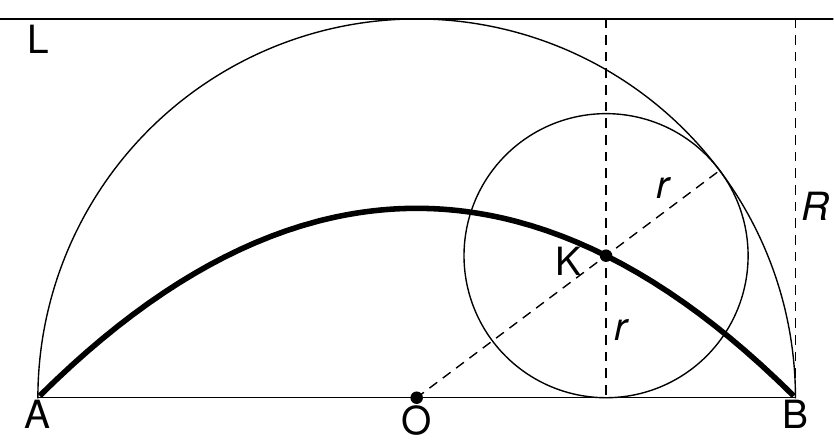}
		\caption{The locus of the centers of circles inscribed in a semicircle.}
		\label{FIG:Locus}
	\end{center}
\end{figure}

Let the semicircle have center~$\mathsf{O}$, radius~{\it\textsf{R}}, and diameter $\mathsf{AB}$, and let an inscribed circle have center~$\mathsf{K}$ and radius~{\it\textsf{r}}. Drawing a line~$\mathsf{L}$ parallel to $\mathsf{AB}$ at a distance~$R$ as in Figure~\ref{FIG:Locus},
we see that
$$
{\rm distance}(\mathsf{K},\mathsf{O}) = \text{{\it\textsf{R}}}- \text{{\it\textsf{r}}} \,= {\rm distance}(\mathsf{K},\mathsf{L}).
$$
Hence $\mathsf{K}$ lies on the parabola with focus $\mathsf{O}$ and directrix $\mathsf{L}$. Since $\mathsf{O}$ lies on $\mathsf{AB}$, it is the latus rectum, implying our claim. Finally, applying it to each semicircle of the arbelos leads to Property~\ref{PROP:locus}.
\end{proof}

\paragraph{Acknowledgments.} I thank both referees for suggesting changes that improved the exposition, and Harold Boas for correcting two numerical errors.

\bigskip

\noindent\textit{209 West 97th Street, New York, NY  10025\\
jsondow@alumni.princeton.edu}

\end{document}